\documentclass{article}
\usepackage[utf8]{inputenc}
\usepackage{amsfonts}
\usepackage{amssymb,amsthm}
\usepackage[margin=1in]{geometry}
\usepackage{amsmath}
\usepackage{fixmath}
\usepackage{hyperref}
\newtheorem{thm}{Theorem}[section]
\newtheorem{rem}[thm]{Remark}
\newtheorem{lem}[thm]{Lemma}
\newtheorem{cor}[thm]{Corollary}
\newtheorem{Def}[thm]{Definition}
\newtheorem{prop}[thm]{Proposition}

\usepackage{authblk}

\hypersetup{hidelinks}

\title{Ramsey--Cass--Koopmans and a Central Bank:\\ An Exploration of Asymptotic Stability}
\author{Sergey Patsuk}
\author{Tyrus Berry}
\affil{Department of Mathematics, George Mason University\\
\href{mailto:spatsuk@gmu.edu}{\texttt{spatsuk@gmu.edu}}, \href{mailto:tberry@gmu.edu}{\texttt{tberry@gmu.edu}}}
\date{November 2026}
\begin{document}
\maketitle
\begin{abstract}
Ramsey-Cass-Koopmans is a foundational neoclassical economic framework that provides a saddle-path stable macro-economic model. We propose a second optimizing agent, a Central Bank that seeks to optimize a public welfare function. We show that this is a controllable system and provide conditions and interpretations for local asymptotic stability. We then interpret the role of economic and business cycles that arise from complex eigenvalues.
\end{abstract}




\section{Introduction}
In 1965, David Cass~\cite{Cass1965} and Tjalling Koopmans~\cite{Koopmans1965} demonstrated that economic variables, such as consumption, could be made endogenous by assuming that they optimize a utility function.  Solving the optimization problem yields a general form for the equations that govern the evolution of these economic variables.  The Ramsey-Cass-Koopmans (RCK) model was the first success of this methodology, which demonstrated an equilibrium that optimally balances capital reinvestment and consumption.  However, the nature of the optimization problem implies that this equilibrium can never be stable.  Instead, under reasonable assumptions on the utility functions, one can show that there exists a so-called \emph{saddle-path} which leads to the unstable equilibrium.  The saddle-path is a mathematical abstraction; an infinitely thin curve which requires the optimizing agent to continuously update their consumption with infinite precision.  Moreover, the saddle-path is itself unstable, so that any perturbation leads to exponentially fast deviation, meaning that even numerical simulations cannot stay on this path.  

A necessary condition for optimality, called the \emph{transversality condition}, shows that the saddle-path is the only globally optimal solution~\cite{Michel1982,Kamihigashi2001}.  This indicates that global optimality requires perfect information and exact mathematical precision, and is thus impossible for any real agent.  In contrast, if we drop the transversality condition, we still get the same equations (e.g. RCK) representing an agent making locally optimal choices.  The only consequence of dropping the transversality condition is that we have to accept that the optimal equilibrium is unstable.  While several economists like Frank Hahn~\cite{Hahn1960,Hahn1966} and Mordecai Kurz~\cite{Kurz1968} still explored these mathematical assumptions in the 1960s and 1970s, this exploration dwindled as these conditions became more or less accepted. However, in the time since, these assumptions have been more and more practically challenging to defend. With technological advances came a greater range of economic conditions, and as the economy became more global, interconnected, and subject to large asymmetric but infrequent shocks, excluding those very paths by transversality and relying on firms and much less individuals for having perfect foresight in economic planning rather than myopic optimization has been challenged~\cite{LeiNoussair2002,CoibionGorodnichenko2015,Manski2004}.

Here, we show that by introducing a second optimizing agent, modeling a central bank, the optimal equilibrium of the RCK system can be stabilized without the transversality condition.  This represents a realistically achievable equilibrium, without requiring any agents to have perfect knowledge or infinite precision.  Specifically, we propose an open-loop Nash type model starting from the original RCK system and adding in a Central Bank which exerts an intervention force on capital output.  Following the methodology of Cass and Koopmans, we make this intervention endogenous by assuming that the Central Bank chooses the force to optimize its own welfare function.  In general the welfare function can be any function that depends on current household capital, consumption, and on the current intervention force. We show that this new system is controllable and, under explicit conditions on the Central Bank's welfare function and household preferences, the optimal equilibrium is locally asymptotically stable (LAS); when $\rho>n$, this eliminates the need to impose the transversality condition (Section~\ref{sec:transversality}). Furthermore, this external agent removes the need for perfect foresight, since under these stability conditions the Central Bank's response corrects deviations from the equilibrium from any nearby initial condition, so that neither households nor the Central Bank need foresee the equilibrium path in order to reach it.  As a stable equilibrium, it is thus robust to noise and reasonable perturbations of the model as well as imprecision in the actions of either optimizing agent.

\section{Economic Background and Intuition}
Through the 20th century a plethora of macro-economic models were created exploring the optimization of consumption given a desire to both save capital as well as pursue production of a good(s). This topic was first written about by Frank Ramsey in 1928~\cite{Ramsey1928}, heavily under the influence of Keynes: exploring optimization of enjoyment, and decoupling that process from capital production. This came at the cost of assuming the exclusion of shocks such as warfare or natural disasters, and as Ramsey put it, an excessively high rate of savings needed to attain maximal utility or "bliss". During the 1940s, a new, "post-Keynesian" approach was born in macro-economic thought which led to what we now know as the Harrod-Domar model. Proposed in two papers, Roy Harrod (1939)~\cite{Harrod1939} and Evsey Domar (1946)~\cite{Domar1946}, this model moved from the classic Keynesian short-run objective into a longer time period; importantly this also moved from a static set of equations, to a dynamical system. Such a contribution allowed for more realistic modeling of the economic reality than the time-adjusted steps in the static system which presupposed an indefinite continuation of observed patterns; which maintains a self-correcting system returning to equilibrium after a capital shock. The dynamical system proposed moved past this, exploring the conditions for a stable, natural rate of growth; but at the cost of self-aggravation, an economic system continuing in the direction of a shock in the period afterward instead of correcting back to equilibrium. This unique path produced by staying at the equilibrium without shocks would later be known and described by Robert Solow and others as a "knife-edge" problem~\cite{Solow1956}. Following this train of thought, Robert Solow (1956)~\cite{Solow1956}, inspired by Neo-classical economics, made a massive break-through, taking the Harrod-Domar system, and augmenting it with an exogenous growth process outside of capital-accumulation, and showing that this allows the economy to naturally tend towards a steady state thus removing the knife-edge. A key addition in the model was a flexible relation of capital to labor, $r=\frac{K}{L}$. Specifically, he notes that this mechanism can fix the downward self-correcting mechanism from Harrod-Domar. An interesting consideration from this model that was later shelved by Cass and Koopmans was that of rentals and interest. Given the possibility of rentals a closed form equation was developed that gave the owner of capital the opportunity at any time to loan out capital for a rate of interest, or convert it into output.  A troubling pattern did emerge from this, as Solow himself admitted, that a system can be very productive, leading to high savings and growth rate, causing $r$ to increase to infinity. On the converse, a system could be very un-productive leading to an $r$ diminishing towards $0$. Finally we get to the contribution of David Cass (1965)~\cite{Cass1965} and Tjalling Koopmans (1965)~\cite{Koopmans1965}, which will be the main focus of our paper. This paper expanded on the Solow assumption of a constant savings rate, by modeling a representative household choice between contributing a unit of capital to output or saving it, discounting future utility at a rate $\rho$, with the intent of modeling utility at an infinite time horizon.
\subsection{On the Stability of Growth Equilibrium}
The question of whether a growth equilibrium is stable did not begin with the optimal growth models, and it did not end with them. Solow (1956)~\cite{Solow1956} removed the knife-edge of Harrod-Domar by allowing the ratio of capital to labor to adjust, but did so in a world with a single good; and it was this assumption, rather than the flexibility of factor proportions alone, that Frank Hahn (1966)~\cite{Hahn1966} would later identify as being of equal importance. Hahn (1960)~\cite{Hahn1960} gave the question its formal shape, defining stability as every root of the linearized system having a negative real part, and showing that whether such an equilibrium is stable depends heavily on how quickly wages and the rate of interest adjust; importantly, one of the cases considered has the interest rate moved by a rule responding to the state of the economy, which makes stability depend on an agent acting on the system from outside. This was done at the cost of assuming savings behavior, rather than deriving it from the choices of households.

A second line of work questioned the process by which an economy reaches equilibrium at all. The usual device, tâtonnement, forbids trade until prices have cleared the market. Hahn and Negishi (1962)~\cite{HahnNegishi1962} dropped this, and showed that when trade is permitted at prices that are not equilibrium prices, every limit point of the adjustment path is an equilibrium, whatever the form of the excess demand functions; this was shown for a pure exchange economy, and so says nothing about an economy that accumulates capital.

The difficulty returned once capital was allowed to take many forms. Hahn (1966)~\cite{Hahn1966}, in a model with one consumption good, several capital goods, and workers who do not save, found that equilibrium paths with correct expectations need not be unique even when initial conditions and expectations are given, and that with two capital goods there is a wide class of initial conditions from which equilibrium paths do not approach balanced growth. He traced this to two features: neoclassical models contain no investment function determining how savings are allocated across capital goods, and the requirement that expectations be both correct and identical across all agents is, in his words, a hard postulate to swallow. Kurz (1968)~\cite{Kurz1968} then generalized the result to discounted optimization, showing that competitive growth processes which are always in equilibrium, and which could be constructed from the optimization of an intertemporally independent utility functional with a constant discount factor, are locally unstable; in the planning interpretation, an error in the initial price system is magnified unless it is corrected. Both authors indicate where a remedy might lie: Hahn in an investment function or a government enforcing an efficient accumulation plan, and Kurz in an adjustment mechanism between spot and futures prices.

The principal positive response is Heller (1975)~\cite{Heller1975}, who showed that in an economy with infinite-horizon Ramsey households, and a market open today for delivery at every future date, price adjustment by an auctioneer according to excess demand is stable in the sense of Lyapunov, in contrast to the saddle-point instability of myopic foresight dynamics. The result is a qualified one: stability is asymptotic only under strong gross substitutability; households are required to discount the future faster than the population grows; and, as Heller himself notes, the auctioneer must know the subjective discount rate, which sidesteps the question of what mechanism re-aims an economy headed for long-run disequilibrium. The space in which the result holds is also infinite-dimensional, which is why Lyapunov and asymptotic stability need not coincide.

Taken together, each of these works recovers stability, or locates its absence, by closing the model in some manner external to the households and firms: a rule for the interest rate, an auctioneer, a full set of futures markets, or a planner who selects the path outright. What remains unexamined, to our knowledge, is a closing agent that is neither a rule nor a planner, but an optimizer in its own right acting on the same capital stock as optimizing households; this is the system we study in the following section.

Brock and Scheinkman (1976)~\cite{BrockScheinkman1976} developed results regarding convergence of optimal growth economic models, exploring necessary conditions for Global and Local Asymptotic stability in Hamiltonian optimal control systems. In particular, they questioned whether the convergence properties of the one-good Cass-Koopmans model could be preserved under modifications of the underlying economic system. Benhabib and Nishimura (1979)~\cite{BenhabibNishimura1979} subsequently demonstrated that optimal growth models with multiple capital goods can admit closed orbits through Hopf bifurcations, showing that intertemporal optimization need not imply convergence to a stationary equilibrium.

These explore stability in optimal growth models under changing structures in regards to long-term behavior. Rather than introducing additional capital goods or altering the original production technology, we consider the introduction of a second optimizing agent, representing a Central Bank, which acts directly on the capital accumulation equation. We then examine whether the resulting coupled system admits locally asymptotically stable equilibria and how changes in the Central Bank's objective influence the stability and potential oscillatory behavior of the system.

\section{A Modified Ramsey--Cass--Koopmans System}
The model developed by Ramsey, Cass, and Koopmans proposed a system solving the following problem:
\[ \max_c \int_0^\infty e^{-(\rho-n)t} u(c) \, dt \]
subject to the ODE constraint,
\[ \dot k = \tilde{f}(k) - (n+\delta)k - c \]
where $k$ is the state variable and $c$ is the control parameter; where $\rho$ is the discount rate, $c$ and $u(c)$ are the rate of consumption and the utility function respectively, $k$ and $\tilde{f}(k)$ are capital and the production function respectively, and finally $n$ and $\delta$ are the population growth rate and the rate of depreciation respectively. An important contribution was that this system is saddle-path stable: the Jacobian at the interior equilibrium has real eigenvalues of opposite sign, so that the set of initial conditions from which the system converges to the steady state forms a one-dimensional stable manifold, and the initial consumption level must be selected so as to lie on it. In the time since its publication several of their assumptions have been questioned.

The first is that the household-firm system is self-contained: no agent outside it alters the aggregate resource constraint, so that capital evolves only through the saving decisions of the households themselves. There is no mechanism by which resources are injected into or withdrawn from the accumulation equation. This assumption is not intrinsic to the optimizing framework. Barro \cite{Barro1990} studies a model in which a government provides services that enter the production technology and are financed by taxation, so that a term reflecting the outside agent's activity appears directly in the aggregate resource constraint, and shows that the resulting growth and welfare properties differ systematically from those of the unaugmented model. Turnovsky \cite{Turnovsky1997} analyzes the transitional dynamics and stability of such a system when the government's flow enters the accumulation equation, which is the question of interest here. We note also that in an infinite-horizon representative-agent economy such interventions may be neutralized by the households themselves, as established by Barro \cite{Barro1974}; the mechanism we propose evades this neutrality by entering the technology rather than the household budget constraint.

The second is that of perfect foresight, which is the assumption that both households and firms have knowledge of the optimal rate of savings and consumption at all times, and that they operate in a rational way, so that the macroeconomic system follows this optimal path with no outside agent enforcing the rule. Evans and Honkapohja \cite{EvansHonkapohja2001} develop the principal alternative, in which agents form expectations by recursive estimation rather than by foresight, and in which expectational stability (E-stability) governs whether a rational expectations equilibrium is attainable as the limit of such a process, and supplies a selection criterion when several equilibria exist. Direct evidence against perfect foresight in this specific model is given by Lei and Noussair \cite{LeiNoussair2002}, who construct an experimental economy with the structure of the Ramsey-Cass-Koopmans model of optimal growth and implement it three ways: decentralized, with multiple agents and a market for capital; with individual subjects placed in the role of social planners; and with a five-agent social planner making a joint decision. Their findings highlight the role of market institutions in facilitating convergence to the optimal steady state, which is to say that convergence is an institutional outcome rather than a property that individual optimizing agents deliver unaided.

We propose a system that relaxes the self-contained system condition, and that is motivated by the empirical failure of perfect foresight, by adding a Central Bank entity that acts on the aggregate resource constraint through capital injections or contractions at a rate $m$. We then show that this system admits an asymptotically stable solution under conditions on the Central Bank's objective and on a single quantity $D$ built from preferences and technology. A Central Bank whose objective ignores consumption can stabilize only when the population growth rate $n$ exceeds the household discount rate $\rho$; one whose objective responds to consumption can do so whenever $D>0$.

\subsection{Economic Intuition}
The Federal Reserve was established in order to furnish an elastic currency and to serve as a lender of last resort; its modern dual mandate of maximum employment and stable prices was formalized later, as the institution's role broadened from a crisis-response facility into an active manager of aggregate demand. Its principal lever is the expansion and contraction of the monetary supply through open market operations in government debt securities, which alter the economy-wide return available on capital and thereby the incentive to save versus consume.

In the reduced-form system we study below, this lever is represented by a single real-side instrument $m$: a rate at which the Central Bank injects capital into, or withdraws it from, the aggregate stock, entering the resource constraint alongside the household's own saving decision. This is a deliberate simplification -- it has no price level, no nominal interest rate distinct from the real return on capital, and no money demand -- but it captures the feature of actual policy most relevant to the stability question posed in the Introduction: a second optimizing agent, with its own objective and its own first-order condition, acting on the same aggregate capital stock the household is trying to accumulate.

We abstract from technological progress, retaining population growth and capital depreciation, and write $f:=\tilde{f}$ for the production function. This allows us to isolate the effects of the Central Bank's intervention on the stability of the capital accumulation system.
\subsection{On Physical and Nominal Capital}
Tobin (1969)~\cite{Tobin1969} described the link between financial and physical capital through $q$, the ratio of the market value of capital to its replacement cost; in long-run equilibrium $q=1$, so that a unit of capital and its nominal claim exchange one for one. This permits two readings of the perturbation $m$: as a Central Bank holding a stock of the single physical capital and buying or selling units of it, or as the Federal Reserve acting on the monetary base, whose effect on investment operates through the deviation of $q$ from $1$. Since our analysis is local to the equilibrium, perturbations move the economy only an arbitrarily small distance, and we require only that $q-1$ be a smooth, invertible function of $m$ near zero; the two readings then differ by a change of variable, which leaves the eigenvalues of the linearized system, and our stability results unchanged.
\section{Interior Equilibrium}
\label{sec:equilibrium}

\subsection{The Dynamical System}
\label{sec:dynamical-system}

Throughout we assume $u'>0$, $u''<0$, and that $f$ is twice continuously differentiable
with $f(0)=0$, $f'>0$, $f''<0$, and $f'$ mapping $(0,\infty)$ onto $(0,\infty)$.

\begin{Def}
\label{Mod-Cass-Koopman}
Let $g$ be a smooth function with $g(0)=0$ and $\gamma:=g'(0)\neq 0$. The system defined by
\[
\max_c \int_0^\infty e^{-(\rho - n)t} u(c)\, dt
\qquad\text{and}\qquad
\max_m \int_0^\infty e^{-(\rho - n) t}\, W(k,c,m)\, dt
\]
\[
\text{s.t.}\qquad \dot k = f(k) + k\,g(m) - (n+\delta)k - c ,
\]
where $W(k,c,m)$ is the objective of the Central Bank, is called the Modified
Ramsey-Cass-Koopmans system. The function $g$ is the \emph{transmission channel}: a
position $m$ of the Central Bank adds $g(m)$ to the return on each unit of capital. It
follows from the original with the addition of the new perturbation parameter $m$, a
reduced-form representation of the public-capital channel of Barro~\cite{Barro1990} and
Turnovsky~\cite{Turnovsky1997} discussed in Section~3.
\end{Def}

\begin{lem}
\label{Equivalence}
The system of Definition~\ref{Mod-Cass-Koopman} is equivalent to the Ramsey-Cass-Koopmans
system of Cass~\cite{Cass1965} when there is no Central Bank intervention.
\end{lem}
\begin{proof}
If $m\equiv 0$ then $g(m)\equiv 0$, and the resource constraint reduces to
$\dot k = f(k)-(n+\delta)k-c$, the original constraint.
\end{proof}

The two maximization problems are not solved by the agents; they are the source of the
rules the agents follow. We make this precise as follows.

\paragraph{Standing assumption.} Each agent satisfies the differential form of its
first-order conditions, but does not solve for its initial level. The initial values
$k(0)$, $c(0)$ and $m(0)$ are inherited, and no transversality condition is imposed on
either agent.

For the household, the Hamiltonian $u(c)+\lambda\bigl(f(k)+k\,g(m)-(n+\delta)k-c\bigr)$
gives $u'(c)=\lambda$ and $\dot\lambda=(\rho-n)\lambda-\lambda\bigl(f'(k)+g(m)-n-\delta\bigr)$;
the population growth rate cancels, and with $\sigma(c):=-u'(c)/u''(c)>0$ the household's
rule is
\begin{equation}
\dot c=\sigma(c)\bigl(f'(k)+g(m)-\delta-\rho\bigr).
\label{eq:household-rule}
\end{equation}
Consumption rises exactly when the return on capital, inclusive of the Central Bank's
position, exceeds the rate of time preference; the rule uses only the current return.

\subsection{Controllability of the Household System}
\label{sec:controllability}

Before specifying how the Central Bank sets $m$, we ask what any instrument entering the
system in this way is able to do. For this subsection only, regard $m(t)$ as a freely
chosen input. With the input held at $m=0$ the pair
\begin{equation}
\dot{k} = f(k)+k\,g(m)-(n+\delta)k-c,
\qquad
\dot{c} = \sigma(c)\bigl(f'(k)+g(m)-\delta-\rho\bigr)
\label{eq:controlled-rck}
\end{equation}
has the rest point $(k^*,c^*)$ of the original system,
\begin{equation}
f'(k^*)=\delta+\rho, \qquad c^*=f(k^*)-(n+\delta)k^* .
\label{eq:equilibrium}
\end{equation}
We write
\[
\kappa:=k^*, \qquad \eta:=-f''(k^*), \qquad \sigma:=\sigma(c^*),
\]
a notation retained for the rest of the paper. With $\Delta x=(k-k^*,\ c-c^*)^{\mathsf T}$,
the linearization of \eqref{eq:controlled-rck} about $(k^*,c^*)$ and $m=0$ is
\begin{equation}
\Delta\dot{x}=M\,\Delta x+b\,m,
\qquad
M=\begin{pmatrix}
\rho-n&-1\\
-\sigma\eta&0
\end{pmatrix},
\qquad
b=\gamma\begin{pmatrix}\kappa\\ \sigma\end{pmatrix}.
\label{eq:linearization}
\end{equation}
Since $\det M=-\sigma\eta<0$, the eigenvalues of $M$ are real and of opposite sign: with
no intervention the rest point is a saddle, and a household that follows
\eqref{eq:household-rule} from an inherited level of consumption does not converge to it.

\begin{thm}[Controllability]
\label{Controllability_theorem}
Let
\begin{equation}
D:=\sigma-\kappa(\rho-n)-\kappa^2\eta .
\label{eq:D}
\end{equation}
The linearized system \eqref{eq:linearization} is controllable if and only if $D\neq 0$.
\end{thm}
\begin{proof}
By the Kalman rank criterion the system is controllable if and only if
$\det[\,b\ \ Mb\,]\neq 0$. A direct computation gives $\det[\,b\ \ Mb\,]=\gamma^2\sigma D$,
and $\gamma\neq0$, $\sigma>0$.
\end{proof}

\begin{rem}
The channel enters Theorem~\ref{Controllability_theorem} only through $\gamma\neq 0$:
whether the instrument can steer the household system is decided by preferences and
technology, through $D$, and not by the shape of $g$. The quantity $D$ compares the
household's willingness to substitute consumption over time, $\sigma$, with the two forces
acting on capital at the rest point, the gap $\rho-n$ and the curvature $\eta$ of the
technology, each weighted by the capital stock.
\end{rem}

\subsection{The Central Bank's Rule}
\label{sec:bank-rule}

The Central Bank's current-value Hamiltonian is
\[
H_F = W(k,c,m) + \lambda_F \bigl( f(k) + k\,g(m) - (n+\delta)k - c \bigr),
\]
in which the path of consumption is taken as given. The first-order condition in $m$ and
the costate equation are
\[
W_m + \lambda_F\,k\,g'(m)=0,
\qquad
\dot\lambda_F=\bigl(\rho+\delta-f'(k)-g(m)\bigr)\lambda_F-W_k .
\]
We eliminate $\lambda_F$ and take $m$ itself as the Central Bank's state: $m$ is
directly observable, whereas $\lambda_F=-W_m/(k\,g'(m))$ is a derived shadow price.

\begin{cor}
\label{Control-system1}
Under the standing assumption, the Modified Ramsey-Cass-Koopmans system is governed by
\begin{align*}
\dot{k} &= f(k) + k\,g(m) - (n+\delta)k - c, \\[0.3em]
\dot{c} &= \sigma(c)\bigl(f'(k) + g(m) - \delta - \rho\bigr), \\[0.3em]
\dot{m} &= \frac{\bigl(\rho+\delta-f'(k)-g(m)\bigr)W_m + k\,g'(m)\,W_k
- \Bigl(W_{mk}-\dfrac{W_m}{k}\Bigr)\dot k - W_{mc}\,\dot c}
{W_{mm}-\dfrac{g''(m)}{g'(m)}\,W_m},
\end{align*}
where $W$ and its derivatives are evaluated at $(k,c,m)$.
\end{cor}
\begin{proof}
Differentiating $\lambda_F\,k\,g'(m)=-W_m$ in time,
\[
\dot\lambda_F\,k\,g'(m)+\lambda_F\,g'(m)\,\dot k+\lambda_F\,k\,g''(m)\,\dot m
=-W_{mk}\dot k-W_{mc}\dot c-W_{mm}\dot m .
\]
Substituting $\lambda_F=-W_m/(k\,g'(m))$ and the costate equation, the left side becomes
\[
-\bigl(\rho+\delta-f'(k)-g(m)\bigr)W_m-k\,g'(m)\,W_k-\frac{W_m}{k}\,\dot k
-\frac{g''(m)}{g'(m)}\,W_m\,\dot m ,
\]
and solving for $\dot m$ gives the stated equation.
\end{proof}

The rule has a direct reading. The Central Bank revises its position in response to the
gap between the rate of time preference and the return on capital, weighted by its
marginal valuation $W_m$ of the instrument; in response to its marginal valuation $W_k$ of
capital, scaled by the strength $k\,g'(m)$ with which the instrument reaches capital; and
in response to the motion of capital and of consumption themselves, through the
cross-partials $W_{mk}$ and $W_{mc}$.

\subsection{Equilibrium Derivation}
\label{sec:equilibrium-sub}

Since neither agent solves for its initial level, an equilibrium of the system is simply
a rest point of Corollary~\ref{Control-system1}.

\begin{lem}
\label{Rest_point_lemma}
Let $k>0$, $g'(m)\neq 0$, and let the denominator of $\dot m$ be nonzero at $(k,c,m)$.
Then $(k,c,m)$ is a rest point if and only if
\[
f'(k)+g(m)=\delta+\rho, \qquad c=f(k)+k\,g(m)-(n+\delta)k, \qquad W_k(k,c,m)=0 .
\]
\end{lem}
\begin{proof}
The first two conditions are $\dot c=0$ and $\dot k=0$. Given these, the numerator of
$\dot m$ reduces to $\bigl(\rho+\delta-f'(k)-g(m)\bigr)W_m+k\,g'(m)\,W_k$; the bracket
vanishes by the first condition, leaving $k\,g'(m)\,W_k=0$, which is $W_k=0$.
\end{proof}

\begin{thm}
\label{Equilibrium_theorem}
Let $k^*$ be the unique solution of $f'(k^*)=\delta+\rho$ and
$c^*:=f(k^*)-(n+\delta)k^*$. If $W_k(k^*,c^*,0)=0$ and $W_{mm}(k^*,c^*,0)\neq 0$, then
$(k^*,c^*,0)$ is a rest point of the system, for every transmission channel $g$. Moreover
$c^*>0$ whenever $\rho\geq n$; when $\rho<n$ we assume in addition that
$f(k^*)>(n+\delta)k^*$.
\end{thm}
\begin{proof}
Existence and uniqueness of $k^*$ follow from the assumptions on $f'$. Since $g(0)=0$,
the triple $(k^*,c^*,0)$ satisfies Lemma~\ref{Rest_point_lemma}. Strict concavity and
$f(0)=0$ give $f(k)>k\,f'(k)$, so that $f(k^*)/k^*>\delta+\rho\geq\delta+n$ when
$\rho\geq n$.
\end{proof}

\begin{rem}
Neither $W_m$ nor $W_c$ appears in Lemma~\ref{Rest_point_lemma}. The first drops out
because at a rest point the return on capital equals the rate of time preference, and it
is this gap that multiplies $W_m$ in the rule. The second appears nowhere in the system:
the Central Bank takes the path of consumption as given, so that only the cross-partials
$W_{mc}$ and $W_{kc}$ of its objective can influence the dynamics. The rest point is that
of the original system; the Central Bank changes the motion around it, and not its
location.
\end{rem}

\begin{rem}[Idle Central Bank]
\label{Idle_remark}
For the remainder of the paper we impose in addition $W_m(k^*,c^*,0)=0$. This is not
needed for the rest point to exist. It says that the shadow price
$\lambda_F^*=-W_m^*/(\kappa\gamma)$ vanishes there, so that the Central Bank has no
first-order reason to hold a position, and it is what allows the linearization below to be
written in terms of the second-order structure of $W$ alone. Our analysis is local to
$(k^*,c^*,0)$; Lemma~\ref{Rest_point_lemma} does not exclude further rest points with
$m\neq0$.
\end{rem}

\subsection{Notation}
\label{sec:notation}

The following shorthand is used throughout, with $\ast$ denoting evaluation at
$(k^*,c^*,0)$. In addition to $\sigma$, $\eta$, $\kappa$ and $\gamma$ above,
\[
A := -\frac{W_{mk}^*}{W_{mm}^*}, \qquad
B := -\frac{1}{W_{mm}^*}, \qquad
E := -\frac{W_{mc}^*}{W_{mm}^*},
\]
\[
C := -W_{kk}^* + \frac{(W_{mk}^*)^2}{W_{mm}^*}, \qquad
G := -W_{kc}^* + \frac{W_{mk}^*W_{mc}^*}{W_{mm}^*}, \qquad
\Gamma := -W_{cc}^* + \frac{(W_{mc}^*)^2}{W_{mm}^*},
\]
and, for the channel,
\[
\tilde A:=\gamma A, \qquad \tilde B:=\gamma^2 B, \qquad \tilde E:=\gamma E .
\]
Here $A$ and $E$ measure how the Central Bank's marginal valuation of its instrument
responds to capital and to consumption, each relative to its curvature in $m$; $C$, $G$
and $\Gamma$ describe the curvature of $W$ in $(k,c)$ that remains once $m$ has been set
by the first-order condition.

\begin{lem}
\label{Admissible_lemma}
The Hessian of $W$ at $(k^*,c^*,0)$ is negative definite if and only if
$B>0$, $C>0$ and $C\,\Gamma>G^2$.
\end{lem}
\begin{proof}
Order the variables as $(m,k,c)$. Negative definiteness is equivalent to $W_{mm}^*<0$,
which is $B>0$, together with negative definiteness of the Schur complement of
$W_{mm}^*$, which is $-\bigl(\begin{smallmatrix} C & G\\ G & \Gamma\end{smallmatrix}\bigr)$.
\end{proof}

We call a second-order structure \emph{admissible} when it satisfies
Lemma~\ref{Admissible_lemma}. Since $\Gamma$ does not enter the linearization below and
any $G$ is admitted by taking $\Gamma>G^2/C$, the parameters governing the dynamics are
$(A,E,G)\in\mathbb{R}^3$, $B>0$, $C>0$, and these may be chosen independently. The case
$W=W(k,m)$ is $E=G=0$.

\subsection{Linearization}
\label{sec:linearization}

All of the stability analysis that follows rests on the single computation of this
subsection.

\begin{prop}
\label{Substitution}
The Jacobian of the system of Corollary~\ref{Control-system1} at $(k^*,c^*,0)$ is
\[
J = \begin{pmatrix}
\rho - n & -1 & \kappa\gamma \\[0.4em]
-\sigma\eta & 0 & \sigma\gamma \\[0.4em]
r_1/\gamma & r_2/\gamma & r_3
\end{pmatrix},
\qquad
\begin{aligned}
r_1&:=\kappa\tilde A^2+(\rho-n)\tilde A+\kappa\tilde B C-\sigma\eta\tilde E,\\
r_2&:=\kappa\bigl(\tilde B G+\tilde A\tilde E\bigr)-\tilde A,\\
r_3&:=\sigma\tilde E .
\end{aligned}
\]
In particular $J=S^{-1}J_0\,S$ with $S=\operatorname{diag}(1,1,\gamma)$, where $J_0$ is the
Jacobian of the proportional channel $g(m)=m$ evaluated at $(\tilde A,\tilde B,C,\tilde E,G)$
in place of $(A,B,C,E,G)$.
\end{prop}
\begin{proof}
The first two rows follow by differentiating $\dot k$ and $\dot c$ directly, using
$f'(k^*)=\delta+\rho$, $g(0)=0$ and $g'(0)=\gamma$; the dependence of $\sigma$ on $c$ does
not enter, since its coefficient vanishes at the rest point.

For the third row, write $\dot m=N/Q$ with $N$ the numerator and $Q$ the denominator of
Corollary~\ref{Control-system1}. At the rest point $N^*=0$ and $Q^*=W_{mm}^*$, so that
$\partial\dot m/\partial x=(\partial N/\partial x)/W_{mm}^*$ for $x\in\{k,c,m\}$. Each of
the four terms of $N$ carries a factor that vanishes at the rest point, and in the first
term both factors vanish; hence only the derivatives of the vanishing factors survive,
\[
\frac{\partial N}{\partial x}\bigg|_*
=\kappa\gamma\,\frac{\partial W_k}{\partial x}\bigg|_*
-W_{mk}^*\,\frac{\partial \dot k}{\partial x}\bigg|_*
-W_{mc}^*\,\frac{\partial \dot c}{\partial x}\bigg|_* .
\]
Dividing by $W_{mm}^*=-1/B$,
\[
\frac{\partial \dot m}{\partial x}\bigg|_*
=-\kappa\gamma B\,\frac{\partial W_k}{\partial x}\bigg|_*
+A\,\frac{\partial \dot k}{\partial x}\bigg|_*
+E\,\frac{\partial \dot c}{\partial x}\bigg|_* .
\]
Inserting the first two rows and $W_{kk}^*=-C-A^2/B$, $W_{kc}^*=-G-AE/B$,
$W_{km}^*=A/B$ gives, for $x=k,c,m$ respectively,
\[
\kappa\gamma(A^2+BC)+(\rho-n)A-\sigma\eta E,
\qquad
\kappa\gamma(BG+AE)-A,
\qquad
\sigma\gamma E,
\]
which are $r_1/\gamma$, $r_2/\gamma$ and $r_3$. The factorization $J=S^{-1}J_0S$ is read
off from the entries.
\end{proof}

\begin{prop}
\label{Reduced_RH}
The characteristic polynomial of $J$ is $\lambda^3+c_1\lambda^2+c_2\lambda+c_3$ with
\begin{align*}
c_1 &= n-\rho-\sigma\tilde E, \\
c_2 &= \sigma(\tilde A-\eta)-\kappa\tilde A\bigl(\rho-n+\kappa\tilde A\bigr)-\kappa^2\tilde B C
+\sigma(\rho-n+\kappa\eta)\tilde E-\sigma\kappa\bigl(\tilde B G+\tilde A\tilde E\bigr), \\
c_3 &= \sigma\kappa\Bigl[\tilde A(\tilde A-\eta)+\tilde B C
+(\rho-n+\kappa\eta)\bigl(\tilde B G+\tilde A\tilde E\bigr)\Bigr].
\end{align*}
All three eigenvalues have negative real part, and the rest point is locally
asymptotically stable, if and only if $c_1>0$, $c_3>0$ and $c_1c_2-c_3>0$.
\end{prop}
\begin{proof}
By Proposition~\ref{Substitution} the matrices $J$ and $J_0$ are similar, so they share a
characteristic polynomial, whose coefficients are $c_1=-\operatorname{tr}J_0$, $c_2$ the
sum of the principal $2\times2$ minors of $J_0$, and $c_3=-\det J_0$:
\[
c_1=n-\rho-r_3,\qquad
c_2=-\kappa r_1-\sigma r_2+(\rho-n)r_3-\sigma\eta,\qquad
c_3=\sigma\bigl[r_1+(\rho-n+\kappa\eta)r_2+\eta r_3\bigr].
\]
Substituting $r_1,r_2,r_3$ gives the stated expressions. The last assertion is the
Routh--Hurwitz criterion for a cubic.
\end{proof}

Three consequences are immediate, and replace separate computations for each channel.

\begin{cor}
\label{Channel_cor}
\begin{enumerate}
\item[(a)] Local stability depends on the transmission channel only through its slope
$\gamma=g'(0)$; no higher derivative of $g$ enters $J$.
\item[(b)] The slope enters only through $\tilde A$, $\tilde B$, $\tilde E$. A system with
channel $g$ and Central Bank $(A,B,C,E,G)$ therefore has the same local dynamics as the
system with the proportional channel and Central Bank $(\gamma A,\gamma^2B,C,\gamma E,G)$.
\item[(c)] The trace is
\[
\operatorname{tr}J=\rho-n+\sigma\gamma E .
\]
When the Central Bank's objective does not depend on consumption, $E=0$, the trace is
$\rho-n$ for every channel and every admissible structure, and no such Central Bank
renders the rest point locally asymptotically stable when $\rho>n$.
\end{enumerate}
\end{cor}
\begin{proof}
Parts (a) and (b) restate Proposition~\ref{Substitution}; (b) is the change of variable
$m\mapsto\gamma m$. For (c), $\operatorname{tr}J=-c_1$, and if $\rho-n>0$ with $E=0$ the
real parts of the eigenvalues sum to a positive number.
\end{proof}

\section{Stability}
\label{sec:stability-impossibility}

By Corollary~\ref{Channel_cor}(b), the local dynamics for a channel of slope $\gamma$
coincide with those of the proportional channel under
$(A,B,E)\mapsto(\tilde A,\tilde B,\tilde E)$. Throughout this section and Section~\ref{sec:hopf} we
therefore set $\gamma=1$ and drop the tildes; every statement holds for general $g$ on
restoring them. We write
\[
a:=\rho-n,
\]
so that the first row of $J$ is $(a,\,-1,\,\kappa)$ and $D=\sigma-\kappa a-\kappa^2\eta$.

\subsection{Realizing the Characteristic Polynomial}

The coefficients $c_1,c_2,c_3$ of Proposition~\ref{Reduced_RH} are affine in the third
row $(r_1,r_2,r_3)$ of $J$, and the third row is the only part of $J$ the Central Bank's
objective can affect. The following lemma converts questions about stability into
questions about which third rows an admissible objective can produce.

\begin{lem}
\label{Realization_lemma}
Let $D\neq 0$ and let $(c_1,c_2,c_3)\in\mathbb{R}^3$. Define
\[
r_3:=-a-c_1, \qquad p:=c_2-a\,r_3+\sigma\eta, \qquad q:=\frac{c_3}{\sigma}-\eta\,r_3,
\]
\[
r_1:=\frac{\sigma q+(a+\kappa\eta)\,p}{D}, \qquad r_2:=-\frac{p+\kappa q}{D}.
\]
Then $J$ has characteristic polynomial $\lambda^3+c_1\lambda^2+c_2\lambda+c_3$ if and
only if its third row is $(r_1,r_2,r_3)$, and this third row is produced by some
admissible second-order structure if and only if
\begin{equation}
r_1+\eta\,r_3>-\frac{a^2}{4\kappa}.
\label{eq:admissible-row}
\end{equation}
\end{lem}

\begin{proof}
From the proof of Proposition~\ref{Reduced_RH},
\[
c_1=-a-r_3, \qquad
p=-\kappa r_1-\sigma r_2, \qquad
q=r_1+(a+\kappa\eta)\,r_2 .
\]
The first fixes $r_3$. The remaining two form a linear system in $(r_1,r_2)$ with
determinant $\sigma-\kappa(a+\kappa\eta)=D\neq0$, whose unique solution is as stated.

For admissibility, Proposition~\ref{Substitution} gives $r_3=\sigma E$,
$r_2=\kappa(BG+AE)-A$ and $r_1=\kappa A^2+aA+\kappa BC-\eta r_3$. The first determines
$E$. Given $E$ and any $A\in\mathbb{R}$, $B>0$, the third determines
$\kappa BC=r_1+\eta r_3-\kappa A^2-aA$ and the second determines $G$; $\Gamma$ is then
chosen with $\Gamma>G^2/C$. Hence the row is admissible if and only if
$r_1+\eta r_3-\kappa A^2-aA>0$ for some $A$. The maximum of $-\kappa A^2-aA$ is
$a^2/(4\kappa)$, attained at $A=-a/(2\kappa)$, which gives \eqref{eq:admissible-row}.
\end{proof}

\subsection{The Trace and the Impossibility Theorem}

\begin{thm}
\label{Impossibility_theorem}
The trace of $J$ is
\[
\operatorname{tr}J=\rho-n+\sigma E .
\]
Consequently:
\begin{enumerate}
\item[(a)] If the Central Bank's objective does not depend on consumption ($E=G=0$) and
$\rho\geq n$, the rest point is not locally asymptotically stable, for every admissible
$(A,B,C)$ and every transmission channel.
\item[(b)] If $\rho\geq n$, local asymptotic stability requires $E<0$, that is
$W_{mc}^*<0$.
\end{enumerate}
\end{thm}

\begin{proof}
The trace is $-c_1$ from Proposition~\ref{Reduced_RH}. Local asymptotic stability
requires $c_1>0$, i.e. $\sigma E<-(\rho-n)$. For (a), $E=0$ makes this $\rho<n$. For (b),
$\rho\geq n$ forces $\sigma E<0$, and since $W_{mm}^*<0$ the sign of $E$ is that of
$W_{mc}^*$.
\end{proof}

\begin{rem}
\label{Additive_remark}
For the additive channel $\dot k=f(k)+g(m)-(n+\delta)k-c$, which lies outside the family
of Definition~\ref{Mod-Cass-Koopman}, the instrument does not enter the household's rule,
the $(2,3)$ and $(3,3)$ entries of the Jacobian vanish, and $\operatorname{tr}J=\rho-n$
whatever the value of $E$. Part (b) of Theorem~\ref{Impossibility_theorem} therefore
cannot be satisfied by a Central Bank that adds resources without acting on the return.
\end{rem}

\begin{rem}
\label{rhoF_remark}
If the Central Bank discounts at $\rho_F-n$ in place of $\rho-n$, the first term of the
numerator of $\dot m$ becomes $(\rho_F+\delta-f'(k)-g(m))W_m$, which no longer vanishes at
the rest point; $(k^*,c^*,0)$ remains a rest point under the idle normalization
$W_m^*=W_k^*=0$ of Remark~\ref{Idle_remark}. The third row of $J$ acquires
$(\rho_F-\rho)(-A,\,-E,\,1)$, and
\[
\operatorname{tr}J=\rho_F-n+\sigma E .
\]
A Central Bank more patient than the population grows, $\rho_F<n$, thus relaxes the
trace condition independently of $E$.
\end{rem}

\subsection{Existence of a Stabilizing Objective}

\begin{thm}
\label{Possibility_theorem}
Suppose $\rho\geq n$. There exists an admissible second-order structure for which the
rest point is locally asymptotically stable if and only if
\[
D=\sigma-\kappa(\rho-n)-\kappa^2\eta>0 .
\]
\end{thm}

\begin{proof}
\emph{Necessity.} Suppose the rest point is locally asymptotically stable, so that
$c_1,c_2,c_3>0$, and note $a\geq0$. With $r_3=-a-c_1$ as in
Lemma~\ref{Realization_lemma},
\[
p=c_2+a(a+c_1)+\sigma\eta>0, \qquad q=\frac{c_3}{\sigma}+\eta(a+c_1)>0 .
\]
From the proof of Lemma~\ref{Realization_lemma}, $p+\kappa q=-D\,r_2$ and
$\sigma q+(a+\kappa\eta)p=D\,r_1$. If $D=0$ the first gives $p+\kappa q=0$, which is
impossible. If $D<0$, then $r_1=-N/|D|$ with $N:=\sigma q+(a+\kappa\eta)p>0$. Now
$|D|=\kappa(a+\kappa\eta)-\sigma<\kappa(a+\kappa\eta)$ and
$N>(a+\kappa\eta)\,a(a+c_1)\geq(a+\kappa\eta)\,a^2$, so
\[
r_1+\eta r_3=-\frac{N}{|D|}-\eta(a+c_1)<-\frac{a^2}{\kappa}\leq-\frac{a^2}{4\kappa}
\]
when $a>0$, while for $a=0$ the left side is negative and the right side is zero. Either
way \eqref{eq:admissible-row} fails, and by Lemma~\ref{Realization_lemma} no admissible
structure produces these coefficients.

\emph{Sufficiency.} Let $D>0$. Fix $c_1>\max\{0,\,-2(a+\kappa\eta)\}$ and, for $t>0$,
take $c_3=t$ and $c_2=2t/c_1$, so that $c_1>0$, $c_3>0$ and $c_1c_2-c_3=t>0$. With
$r_3=-a-c_1$ fixed, Lemma~\ref{Realization_lemma} gives
\[
r_1=\frac{t}{D}\Bigl(1+\frac{2(a+\kappa\eta)}{c_1}\Bigr)+\text{const},
\]
where the constant does not depend on $t$. The bracket is positive by the choice of
$c_1$, so $r_1\to+\infty$ as $t\to\infty$, and \eqref{eq:admissible-row} holds for $t$
large. The resulting structure is admissible and the Routh--Hurwitz conditions hold.
\end{proof}

\begin{rem}
The sufficiency argument does not use $\rho\geq n$: $D>0$ is sufficient in every regime.
\end{rem}

When the Central Bank's objective does not depend on consumption, the trace confines
stability to $\rho<n$, and there the condition sharpens.

\begin{thm}
\label{Possibility_kl}
Let $E=G=0$ and write $d:=n-\rho$. There exists an admissible $(A,B,C)$ for which the
rest point is locally asymptotically stable if and only if
\[
d>\kappa\eta, \qquad\text{i.e.}\qquad n-\rho>-k^*f''(k^*).
\]
Every stabilizing structure has
\[
\frac{\sigma\eta}{\sigma+\kappa d}<A<\frac{d}{\kappa},
\]
so in particular $A>0$. Conversely, if $d>\kappa\eta$, the choice $A=\eta$, any $B>0$,
and
\[
0<C<C^*:=\frac{\eta\,d\,(d-\kappa\eta)}{B(\kappa d+\sigma)}
\]
is stabilizing.
\end{thm}

\begin{proof}
With $E=G=0$, Proposition~\ref{Reduced_RH} gives $c_1=d$,
$c_3=\sigma\kappa\bigl[A(A-\eta)+BC\bigr]$, and, writing $s:=\sigma+\kappa d$,
\begin{equation}
c_1c_2-c_3=(d-\kappa A)(sA-\sigma\eta)-\kappa s\,BC .
\label{eq:hurwitz-kl}
\end{equation}
Stability requires $d>0$ and $(d-\kappa A)(sA-\sigma\eta)>\kappa sBC>0$.

If $\kappa A>d$, this requires $sA<\sigma\eta$, hence $A<\eta$; but then
$c_2=\sigma(A-\eta)+\kappa A(d-\kappa A)-\kappa^2BC<0$, a contradiction. If $\kappa A=d$
the left side of \eqref{eq:hurwitz-kl} is negative. Hence $\kappa A<d$ and
$sA>\sigma\eta$, which is the stated range of $A$.

Suppose now $d\leq\kappa\eta$. Then $A<d/\kappa\leq\eta$, and $c_3>0$ requires
$BC>A(\eta-A)$. Substituting into \eqref{eq:hurwitz-kl},
\[
(d-\kappa A)(sA-\sigma\eta)>\kappa sA(\eta-A),
\]
which simplifies to $d\bigl[A(\sigma+\kappa d-\kappa^2\eta)-\sigma\eta\bigr]>0$. Since
$\kappa d-\kappa^2\eta\leq0$ and $A<\eta$, the bracket is less than $\sigma A-\sigma\eta<0$,
a contradiction. Hence $d>\kappa\eta$ is necessary.

Conversely, let $d>\kappa\eta$ and $A=\eta$. Then $c_3=\sigma\kappa BC>0$ and
\eqref{eq:hurwitz-kl} becomes $\kappa\eta\,d(d-\kappa\eta)-\kappa sBC$, which is positive
exactly when $C<C^*$. Since $c_1,c_3>0$ and $c_1c_2>c_3$, also $c_2>0$, and the
Routh--Hurwitz conditions hold.
\end{proof}

\begin{rem}
Since $d>\kappa\eta$ implies $D=\sigma+\kappa d-\kappa^2\eta>\sigma>0$, the region of
Theorem~\ref{Possibility_kl} lies inside that of Theorem~\ref{Possibility_theorem}.
Allowing the objective to depend on consumption therefore enlarges the stabilizable
region in two directions: it admits $\rho\geq n$, and for $\rho<n$ it admits technologies with
$\kappa\eta\geq d$ but $D>0$, for which an objective in $(k,m)$ alone cannot
stabilize.
\end{rem}
\section{Transversality}
\label{sec:transversality}

In the original Ramsey--Cass--Koopmans system the transversality condition does real
work: together with feasibility, it selects the single initial level of consumption from
which the economy converges. Under the standing assumption of
Section~\ref{sec:dynamical-system} no such selection is made, since $c(0)$ and $m(0)$ are
inherited. We show that on the trajectories the agents actually follow, both agents'
transversality conditions hold automatically when $\rho>n$ and the rest point is locally
asymptotically stable. The conditions remain necessary for optimality, but they no longer
restrict the initial conditions.

Throughout this section, $x(t)=(k(t),c(t),m(t))$ denotes a trajectory of the system of
Corollary~\ref{Control-system1}, $x^*=(k^*,c^*,0)$, and
\[
\lambda_F(x):=-\frac{W_m(k,c,m)}{k\,g'(m)},
\]
which is the Central Bank's shadow price on capital expressed through its first-order
condition. Since $k^*>0$ and $g'(0)=\gamma\neq0$, $\lambda_F$ is smooth near $x^*$, and by
Remark~\ref{Idle_remark}, $\lambda_F(x^*)=0$.

\begin{Def}
\label{TVC_def}
Along a trajectory $x(t)$, the \emph{household} and \emph{Central Bank} transversality
conditions are
\begin{equation}
\lim_{t\to\infty} e^{-(\rho-n)t}\,u'\bigl(c(t)\bigr)\,k(t)=0,
\qquad
\lim_{t\to\infty} e^{-(\rho-n)t}\,\lambda_F\bigl(x(t)\bigr)\,k(t)=0 .
\label{eq:tvc}
\end{equation}
\end{Def}

\subsection{The Household}

\begin{prop}
\label{TVC_household}
Let $x(t)\to x^*$. Then the household transversality condition holds if and only if
$\rho>n$.
\end{prop}

\begin{proof}
By continuity, $u'(c(t))\,k(t)\to u'(c^*)\,k^*$, which is strictly positive since $u'>0$,
$k^*>0$ and $c^*>0$ (Theorem~\ref{Equilibrium_theorem}). Hence
$e^{-(\rho-n)t}u'(c(t))k(t)$ tends to $0$ if $\rho>n$, to $u'(c^*)k^*\neq0$ if $\rho=n$,
and to $+\infty$ if $\rho<n$.
\end{proof}

\subsection{The Central Bank}

The Central Bank's shadow price tends to zero along any convergent trajectory, so its
condition can hold even when $\rho\leq n$; whether it does depends on how fast the
trajectory converges.

\begin{lem}
\label{Decay_lemma}
Let $x^*$ be locally asymptotically stable, and let $\mu<0$ be the largest real part among
the eigenvalues of $J$. For every $\varepsilon>0$ and every trajectory $x(t)$ in the basin
of attraction of $x^*$ there is a constant $K$ such that
\[
|x(t)-x^*|\leq K e^{(\mu+\varepsilon)t}
\qquad\text{and}\qquad
|\lambda_F(x(t))|\leq K e^{(\mu+\varepsilon)t}
\qquad\text{for all }t\geq0 .
\]
\end{lem}

\begin{proof}
Write $y=x-x^*$, so that $\dot y=Jy+R(y)$ with $|R(y)|=o(|y|)$. Since every eigenvalue of
$J$ has real part at most $\mu$, there is an inner product $\langle\cdot,\cdot\rangle$ on
$\mathbb{R}^3$ with $\langle Jy,y\rangle\leq(\mu+\varepsilon/2)|y|^2$. Choose $\delta>0$ so
that $|R(y)|\leq(\varepsilon/2)|y|$ for $|y|<\delta$. On this ball,
$\tfrac{d}{dt}|y|\leq(\mu+\varepsilon)|y|$, so a trajectory that enters it at time $t_0$
satisfies $|y(t)|\leq|y(t_0)|e^{(\mu+\varepsilon)(t-t_0)}$ thereafter. A trajectory in the
basin enters the ball in finite time, and on $[0,t_0]$ it is bounded; this gives the first
bound after enlarging $K$. Since $\lambda_F$ is smooth with $\lambda_F(x^*)=0$, it is
Lipschitz near $x^*$, which gives the second.
\end{proof}

\begin{prop}
\label{TVC_bank}
Let $x^*$ be locally asymptotically stable and let $\mu$ be as in
Lemma~\ref{Decay_lemma}. If
\begin{equation}
\mu<\rho-n,
\label{eq:bank-tvc-condition}
\end{equation}
then every trajectory in the basin of attraction of $x^*$ satisfies the Central Bank
transversality condition. In particular this holds whenever $\rho\geq n$.
\end{prop}

\begin{proof}
Choose $\varepsilon>0$ with $\mu+\varepsilon<\rho-n$. By Lemma~\ref{Decay_lemma}, and
since $k(t)$ is bounded along the trajectory,
\[
\bigl|e^{-(\rho-n)t}\lambda_F(x(t))\,k(t)\bigr|
\leq K' e^{(\mu+\varepsilon-(\rho-n))t}\to 0 .
\]
When $\rho\geq n$, condition \eqref{eq:bank-tvc-condition} holds because $\mu<0$.
\end{proof}

\begin{cor}
\label{TVC_bank_rho_less_n}
Suppose $\rho<n$. Condition \eqref{eq:bank-tvc-condition} requires
\[
\sigma\gamma E<-2(n-\rho).
\]
In particular, if the Central Bank's objective does not depend on consumption ($E=0$),
condition \eqref{eq:bank-tvc-condition} fails at every locally asymptotically stable rest
point.
\end{cor}

\begin{proof}
The real parts of the eigenvalues sum to $\operatorname{tr}J=\rho-n+\sigma\gamma E$
(Corollary~\ref{Channel_cor}(c)), so the largest of them is at least one third of
the sum:
\[
\mu\geq\tfrac13\bigl(\rho-n+\sigma\gamma E\bigr).
\]
Condition \eqref{eq:bank-tvc-condition} therefore requires
$\tfrac13(\rho-n+\sigma\gamma E)<\rho-n$, which is $\sigma\gamma E<-2(n-\rho)$. With
$E=0$ this reads $0<-2(n-\rho)$, which is false.
\end{proof}

\begin{rem}
\label{TVC_bank_generic}
Condition \eqref{eq:bank-tvc-condition} is sufficient but not necessary for an
individual trajectory: a trajectory that approaches $x^*$ tangent to a faster eigenspace
can satisfy the Central Bank's condition even when $\mu\geq\rho-n$. Such trajectories are
exceptional. The gradient of $\lambda_F$ at $x^*$ is
$(-A,-E,1)/(\kappa\gamma B)\neq0$, so for a generic trajectory the component of
$\lambda_F(x(t))$ along the slowest eigenspace is nonzero and decays at rate $\mu$; if
$\mu>\rho-n$, the Central Bank's condition then fails.
\end{rem}

\subsection{Both Agents}

\begin{thm}
\label{TVC_theorem}
Let $x^*$ be locally asymptotically stable.
\begin{enumerate}
\item[(a)] If $\rho>n$, every trajectory in the basin of attraction of $x^*$ satisfies
both transversality conditions~\eqref{eq:tvc}.
\item[(b)] If $\rho\leq n$, no trajectory converging to $x^*$ satisfies the household
transversality condition.
\end{enumerate}
By Theorem~\ref{Possibility_theorem}, the hypothesis of (a) can be met by an admissible
Central Bank if and only if $D>0$.
\end{thm}

\begin{proof}
Part (a) combines Propositions~\ref{TVC_household} and~\ref{TVC_bank}; part (b) is
Proposition~\ref{TVC_household}. The final statement is
Theorem~\ref{Possibility_theorem}.
\end{proof}

\begin{cor}
\label{TVC_no_selection}
Suppose $\rho>n$ and $x^*$ is locally asymptotically stable. For every $k(0)$
sufficiently close to $k^*$, the set of inherited initial conditions $(c(0),m(0))$ from
which the economy converges while satisfying both transversality conditions contains an
open neighbourhood of $(c^*,0)$.
\end{cor}

\begin{proof}
The basin of attraction of a locally asymptotically stable rest point is open and
contains $x^*$, so its slice at fixed $k(0)$ near $k^*$ contains an open neighbourhood of
$(c^*,0)$. By Theorem~\ref{TVC_theorem}(a) every trajectory starting there satisfies both
conditions.
\end{proof}

\begin{rem}
\label{TVC_contrast}
Corollary~\ref{TVC_no_selection} is the precise sense in which the transversality
condition is no longer needed. In the original system, $m\equiv0$ and the rest point is a
saddle (Section~\ref{sec:controllability}): for given $k(0)$, convergence occurs from a
single value of $c(0)$. Trajectories starting below it over-accumulate capital and violate
the transversality condition, while those starting above it exhaust the capital stock in
finite time and violate feasibility~\cite{Cass1965}. Since the transversality condition is
necessary for optimality under standard conditions~\cite{Michel1982,Kamihigashi2001}, it
is binding there: it is what pins down $c(0)$. Here it holds on an open set of initial
conditions, and so it selects nothing.
\end{rem}

\begin{rem}
\label{TVC_scope}
When $\rho\leq n$ the household's condition fails on every convergent trajectory
(Theorem~\ref{TVC_theorem}(b)). The convergent paths of
Theorem~\ref{Possibility_kl}, which require $\rho<n$, are therefore paths of agents
following their rules, but not optimal paths for the household. The interpretation of the
system as the dynamics of two optimizing agents is available only when $\rho>n$.
\end{rem}
\section{Bifurcations}
\label{sec:hopf}

\subsection{The Boundary of the Stable Region}

\begin{prop}
\label{No_crossing}
Along any continuous path in parameter space that starts in the stable region and on
which $c_3>0$, the quantity $c_1c_2-c_3$ reaches zero no later than $c_1$ does.
Consequently the boundary of the stable region consists of the surface $c_3=0$, where a
real eigenvalue passes through zero, and the surface $c_1c_2=c_3$ with $c_1>0$; the
surface $c_1=0$ meets the boundary only where $c_3=0$.
\end{prop}

\begin{proof}
At any point with $c_1=0$ and $c_3>0$, $c_1c_2-c_3=-c_3<0$. Since $c_1c_2-c_3>0$ at the
start of the path, it vanishes before $c_1$ does. An eigenvalue $\lambda=i\omega$ satisfies
$c_3-c_1\omega^2=0$ and $\omega(c_2-\omega^2)=0$; for $\omega\neq0$ and $c_1>0$ this is
$c_1c_2=c_3$, and for $\omega=0$ it is $c_3=0$.
\end{proof}

In particular, the trace condition of Theorem~\ref{Impossibility_theorem} is never the
condition through which stability is lost while $c_3>0$: the economy loses stability
through oscillation first.

\subsection{Hopf Bifurcation}

\begin{thm}
\label{Hopf_locus}
Let $c_1>0$. On the surface $c_1c_2=c_3$ with $c_3>0$, the characteristic polynomial
factors as
\begin{equation}
\label{eq:complex}
\lambda^3+c_1\lambda^2+c_2\lambda+c_3=(\lambda+c_1)(\lambda^2+c_2)
\end{equation}
and the spectrum of $J$ is $\{-c_1,\ \pm i\omega_0\}$ with $\omega_0:=\sqrt{c_2}>0$.
\end{thm}

\begin{proof}
Expanding the product reproduces the cubic exactly when $c_3=c_1c_2$. Since $c_1,c_3>0$,
$c_2=c_3/c_1>0$.
\end{proof}

\begin{thm}
\label{Hopf_transversality}
Hold all parameters except $C$ fixed, and let $C=C^*$ be a point of the surface of
Theorem~\ref{Hopf_locus}. The eigenvalue $\lambda(C)$ passing through $i\omega_0$
satisfies
\[
\frac{d}{dC}\operatorname{Re}\lambda(C)\bigg|_{C=C^*}
=\frac{\kappa B\,(\sigma+\kappa c_1)}{2\,(\omega_0^2+c_1^2)}>0,
\]
and the real eigenvalue $\lambda_r(C)$ passing through $-c_1$ satisfies
$\lambda_r'(C^*)=-2\,\frac{d}{dC}\operatorname{Re}\lambda(C)\big|_{C^*}$. The crossing is
transversal: as $C$ increases through $C^*$ the complex pair moves from the left to the
right half-plane while the real eigenvalue remains negative.
\end{thm}

\begin{proof}
The parameter $C$ enters $J$ only through the term $\kappa BC$ in $r_1$. Hence $c_1$ does
not depend on $C$, while $dc_2/dC=-\kappa^2B$ and $dc_3/dC=\sigma\kappa B$. Let
$P(\lambda,C):=\lambda^3+c_1\lambda^2+c_2(C)\lambda+c_3(C)$. By the implicit function
theorem $d\lambda/dC=-P_C/P_\lambda$, and at $(i\omega_0,C^*)$, using $c_2^*=\omega_0^2$,
\[
P_C=\kappa B(\sigma-i\kappa\omega_0), \qquad
P_\lambda=-2\omega_0^2+2ic_1\omega_0=2\omega_0(ic_1-\omega_0).
\]
Multiplying numerator and denominator of $-P_C/P_\lambda$ by $\overline{P_\lambda}$ and
taking the real part gives the stated expression. Since the sum of the eigenvalues is
$-c_1$, independent of $C$, $\lambda_r'+2\operatorname{Re}\lambda'=0$; and $\lambda_r(C^*)=-c_1<0$,
so $\lambda_r$ remains negative near $C^*$.
\end{proof}

\begin{cor}
\label{Hopf_explicit}
In the setting of Theorem~\ref{Possibility_kl} with $A=\eta$, the upper bound $C^*$ is a
Hopf point, at which
\[
\omega_0^2=c_2^*=\frac{\sigma\kappa\eta\,(d-\kappa\eta)}{\kappa d+\sigma}>0 .
\]
For $C$ slightly below $C^*$ the rest point is a stable focus, and slightly above it an
unstable focus, with a real eigenvalue near $-d$.
\end{cor}

\begin{proof}
At $A=\eta$, $c_2=\kappa\eta(d-\kappa\eta)-\kappa^2BC$. Substituting $C^*$,
\[
c_2^*=\kappa\eta(d-\kappa\eta)\Bigl[1-\frac{\kappa d}{\kappa d+\sigma}\Bigr]
=\frac{\sigma\kappa\eta(d-\kappa\eta)}{\kappa d+\sigma},
\]
which is positive since $d>\kappa\eta$. The remaining claims follow from
Theorems~\ref{Hopf_locus} and~\ref{Hopf_transversality}.
\end{proof}

\begin{rem}
By Corollary~\ref{Channel_cor}(b), Theorems~\ref{Hopf_locus}--\ref{Hopf_transversality}
hold for a general channel with $B$ replaced by $\tilde B=\gamma^2B$. The location of the
Hopf point and the direction of crossing depend on $g$ only through $\gamma$; whether the
bifurcation is supercritical or subcritical depends on the second- and third-order terms of
the system, through which $g''(0)$ and $g'''(0)$ enter.
\end{rem}

\subsection{The Zero-Hopf Point}

\begin{prop}
\label{Zero_Hopf}
Suppose $c_1=0$, $c_3=0$ and $c_2>0$. Then the spectrum of $J$ is
$\{0,\pm i\sqrt{c_2}\}$. The two defining conditions are transversal in the parameters
$(\rho,C)$:
\[
\det\frac{\partial(c_1,c_3)}{\partial(\rho,C)}=-\sigma\kappa B\neq0 .
\]
\end{prop}

\begin{proof}
With $c_1=c_3=0$ the characteristic polynomial is $\lambda(\lambda^2+c_2)$. From
Proposition~\ref{Reduced_RH}, $c_1=n-\rho-\sigma E$, so $\partial c_1/\partial\rho=-1$ and
$\partial c_1/\partial C=0$; and $\partial c_3/\partial C=\sigma\kappa B$.
\end{proof}

\begin{rem}
The rest point $(k^*,c^*,0)$ exists for every value of the parameters, so the zero
eigenvalue at $c_3=0$ is not accompanied by a fold of rest points, and the generic
fold-Hopf unfolding does not apply directly. The local dynamics near this point, including
any invariant tori, require the center-manifold normal form, which we do not compute.
\end{rem}
\section{Policy-Induced Cycles}
\label{sec:cycles}
 
We now ask when convergence to the rest point is oscillatory, that is, when $J$ has a
complex pair of eigenvalues while all three eigenvalues have negative real part. Such
damped oscillations are the model's counterpart of business cycles. As in
Section~\ref{sec:stability-impossibility} we set $\gamma=1$ and write $a:=\rho-n$;
every statement holds for a general channel on restoring the tildes.
 
\subsection{The Baseline Has No Cycles}
 
\begin{prop}
\label{No_baseline_cycles}
Without intervention, the linearized household system \eqref{eq:linearization} with
$m\equiv0$ has two real eigenvalues. Hence every oscillation of the linearized
Modified Ramsey--Cass--Koopmans system is induced by the Central Bank.
\end{prop}
 
\begin{proof}
The eigenvalues of $M$ solve $\lambda^2-(\operatorname{tr}M)\lambda+\det M=0$, with
$\operatorname{tr}M=\rho-n$ and $\det M=-\sigma\eta<0$. The discriminant is
$(\rho-n)^2+4\sigma\eta>0$, so the eigenvalues are real and distinct.
\end{proof}
 
\subsection{Risk Aversion and the Existence of Cycles}
 
For the remainder of this subsection let $f(k)=k^\alpha$ with $0<\alpha<1$, and let $u$
have constant relative risk aversion $\theta>0$, so that $\sigma(c)=c/\theta$.
 
\begin{prop}
\label{Theta_bound}
Under these specifications,
\[
\frac{D}{\kappa}
=\frac1\theta\Bigl(\frac{\delta+\rho}{\alpha}-(n+\delta)\Bigr)-(\rho-n)-(1-\alpha)(\delta+\rho).
\]
Consequently, if $\rho\geq n$ and $c^*>0$, then $D>0$ if and only if
\[
\theta<\Theta:=\frac{\dfrac{\delta+\rho}{\alpha}-(n+\delta)}{(\rho-n)+(1-\alpha)(\delta+\rho)} .
\]
\end{prop}
 
\begin{proof}
Since $f'(k)=\alpha k^{\alpha-1}$ and $f''(k)=\alpha(\alpha-1)k^{\alpha-2}$, we have
$k\,f''(k)=(\alpha-1)f'(k)$, and by the rest-point condition $f'(k^*)=\delta+\rho$,
\[
\kappa\eta=-k^*f''(k^*)=(1-\alpha)(\delta+\rho).
\]
Similarly $f(k)/k=k^{\alpha-1}=f'(k)/\alpha$, so that
\[
\frac{c^*}{k^*}=\frac{f(k^*)}{k^*}-(n+\delta)=\frac{\delta+\rho}{\alpha}-(n+\delta).
\]
Since $\sigma=c^*/\theta$, dividing $D=\sigma-\kappa(\rho-n)-\kappa^2\eta$ by $\kappa>0$ gives
the stated expression. When $\rho\geq n$ the denominator of $\Theta$ is positive, and
$c^*>0$ makes the numerator positive, so $D/\kappa>0$ rearranges to $\theta<\Theta$.
\end{proof}
 
\begin{rem}
\label{Theta_remark}
The level of capital drops out: $\Theta$ depends only on the technology share $\alpha$,
the depreciation rate, the rate of time preference and the rate of population growth.
By Theorem~\ref{Possibility_theorem}, when $\rho\geq n$ no admissible Central Bank can
stabilize an economy whose households are more risk-averse than $\Theta$. Because $\Theta$
depends only on $(\alpha,\delta,\rho,n)$, the bound can be compared directly with
evidence on $\theta$. That evidence is far from settled: common calibrations take values
such as $2.5$, $5$ and $10$, while estimates from the consumption Euler equation, once
corrected for publication bias, centre near $1$~\cite{ElminejadHavranekIrsova2022}.
Whether a given economy satisfies $\theta<\Theta$ is therefore sensitive to which of these
is adopted.
\end{rem}
 
\begin{cor}
\label{Cycles_need_theta}
If $\rho\geq n$, the linearized system can exhibit damped oscillations about the rest
point only if $D>0$; under the specifications of Proposition~\ref{Theta_bound}, only if
$\theta<\Theta$.
\end{cor}
 
\begin{proof}
Damped oscillation requires the rest point to be locally asymptotically stable, which
by Theorem~\ref{Possibility_theorem} requires $D>0$.
\end{proof}
 
\subsection{Which Cycles Can an Optimizing Central Bank Produce?}
 
We parametrize stable spectra with a complex pair by a real eigenvalue $-r$ and a pair
$\mu\pm i\omega$, with $r>0$, $\mu<0$ and $\omega>0$. Then
\begin{equation}
c_1=r-2\mu,\qquad c_2=\mu^2+\omega^2-2r\mu,\qquad c_3=r(\mu^2+\omega^2),
\label{eq:cycle-coefficients}
\end{equation}
the period of the oscillation is $T=2\pi/\omega$, and its amplitude is multiplied by
$e^{2\pi\mu/\omega}$ over each period. By Lemma~\ref{Realization_lemma}, such a spectrum
is produced by an admissible Central Bank if and only if
$r_1+\eta r_3>-a^2/(4\kappa)$, with $r_1,r_3$ computed from
\eqref{eq:cycle-coefficients}.
 
\begin{lem}
\label{Fast_cycles}
Suppose $\rho\geq n$ and $D>0$, and fix $r>0$ and $\mu<0$. Then
\[
r_1=\frac{r+a+\kappa\eta}{D}\,\omega^2+R,
\]
where $R$ does not depend on $\omega$, while $r_3$ does not depend on $\omega$.
Consequently there is $\bar\omega$ such that every spectrum $\{-r,\ \mu\pm i\omega\}$
with $\omega>\bar\omega$ is produced by an admissible Central Bank.
\end{lem}
 
\begin{proof}
By \eqref{eq:cycle-coefficients}, $c_1$ does not depend on $\omega$, hence neither does
$r_3=-a-c_1$. In $p=c_2-a\,r_3+\sigma\eta$ the only term depending on $\omega$ is the
$\omega^2$ in $c_2$, and in $q=c_3/\sigma-\eta\,r_3$ it is $r\omega^2/\sigma$. Therefore
\[
\sigma q+(a+\kappa\eta)\,p=(r+a+\kappa\eta)\,\omega^2+\text{const},
\]
and dividing by $D$ gives the expression for $r_1$. Since $r>0$, $a\geq0$ and
$\kappa\eta>0$, the coefficient of $\omega^2$ is positive, so $r_1\to+\infty$ as
$\omega\to\infty$ while $r_3$ is fixed, and the admissibility inequality holds for all
sufficiently large $\omega$.
\end{proof}
 
Lemma~\ref{Fast_cycles} leaves open whether slow cycles are also available. They are:
when $\rho\geq n$ and $D>0$, concavity of the Central Bank's objective places no
restriction on the stable spectra it can produce.
 
\begin{thm}
\label{All_cycles}
Suppose $\rho\geq n$ and $D>0$. Then every cubic $\lambda^3+c_1\lambda^2+c_2\lambda+c_3$
satisfying the Routh--Hurwitz conditions $c_1>0$, $c_3>0$, $c_1c_2>c_3$ is the
characteristic polynomial of $J$ for some admissible second-order structure. In
particular, for every $r>0$, $\mu<0$ and $\omega>0$, the spectrum
$\{-r,\ \mu\pm i\omega\}$ is produced by an admissible Central Bank: damped cycles of
every period and every rate of decay are attainable.
\end{thm}
 
\begin{proof}
Let the coefficients satisfy the Routh--Hurwitz conditions, so that also $c_2>0$, and
define $r_1,r_3,p,q$ as in Lemma~\ref{Realization_lemma}. Since $r_3=-(a+c_1)$,
\[
p=c_2+a(a+c_1)+\sigma\eta, \qquad \sigma q=c_3+\sigma\eta(a+c_1),
\]
and $\eta r_3=-\eta(a+c_1)$. Writing $r_1+\eta r_3$ over the common denominator $D$,
its numerator is
\[
\sigma q+(a+\kappa\eta)\,p-D\,\eta(a+c_1).
\]
Using $D=\sigma-\kappa(a+\kappa\eta)$, the terms $\sigma\eta(a+c_1)$ cancel, and the
numerator becomes
\[
c_3+(a+\kappa\eta)(c_2+\sigma\eta)+(a+c_1)\bigl[a(a+\kappa\eta)+\kappa\eta(a+\kappa\eta)\bigr]
=c_3+(a+\kappa\eta)(c_2+\sigma\eta)+(a+\kappa\eta)^2(a+c_1).
\]
Since $a\geq0$, every factor is positive, so
\[
r_1+\eta r_3=\frac{c_3+(a+\kappa\eta)(c_2+\sigma\eta)+(a+\kappa\eta)^2(a+c_1)}{D}>0
\geq-\frac{a^2}{4\kappa},
\]
and by Lemma~\ref{Realization_lemma} the coefficients are produced by an admissible
structure. The final assertion applies this to \eqref{eq:cycle-coefficients}, which
satisfy the Routh--Hurwitz conditions because all roots have negative real part.
\end{proof}
 
\begin{cor}
\label{Cycles_iff}
Suppose $\rho\geq n$. Damped oscillations about the rest point can be induced by an
admissible Central Bank if and only if $D>0$; under the specifications of
Proposition~\ref{Theta_bound}, if and only if $\theta<\Theta$. When they can, every
period $T=2\pi/\omega$ and every rate of decay is attainable.
\end{cor}
 
\begin{proof}
Necessity is Corollary~\ref{Cycles_need_theta}; sufficiency and the final statement are
Theorem~\ref{All_cycles}.
\end{proof}
 
\begin{rem}
\label{Cycles_interpretation}
Two consequences are worth stating. First, the length of the cycle is not determined by
preferences or technology: once $D>0$, it is set entirely by the second-order structure
of the Central Bank's objective. The model therefore does not predict a business-cycle
frequency; it predicts that cycles exist only through policy, and only when households
are not too risk-averse. Second, Theorem~\ref{All_cycles} contains Lemma~\ref{Fast_cycles}
and yields a direct proof of the sufficiency part of Theorem~\ref{Possibility_theorem}
when $\rho\geq n$: any Hurwitz cubic, oscillatory or not, is attainable.
\end{rem}
 
\begin{rem}
\label{Cycles_Hopf_remark}
By Proposition~\ref{No_crossing}, a stable rest point loses stability, while $c_3>0$,
only through the surface $c_1c_2=c_3$, on which $J$ has the purely imaginary pair
$\pm i\omega_0$ of Theorem~\ref{Hopf_locus}. Near that surface the rest point is a
stable focus, so that damped cycles of period close to $2\pi/\omega_0$ are the generic
form taken by the dynamics just before stability is lost.
\end{rem}
\section{Comparison with Unrestricted Feedback Rules}
\label{sec:feedback}

The Central Bank's rule of Corollary~\ref{Control-system1} is restricted by the concavity
of its objective. We compare it with a rule subject to no such restriction: a linear
feedback in the state of the household system \eqref{eq:linearization},
\begin{equation}
m=-K_1(k-k^*)-K_2(c-c^*).
\label{eq:feedback}
\end{equation}
The closed-loop matrix is $M-bK$, and with $\hat K_i:=\gamma K_i$,
\[
\operatorname{tr}(M-bK)=a-\kappa\hat K_1-\sigma\hat K_2, \qquad
\det(M-bK)=-\sigma\bigl[\eta+\hat K_1+(a+\kappa\eta)\hat K_2\bigr].
\]

\begin{thm}[Rule existence]
\label{Rule_theorem}
There exist gains $(K_1,K_2)$ for which \eqref{eq:feedback} makes $M-bK$ Hurwitz if and
only if $D\neq0$. For given $\hat K_2$, the stabilizing values of $\hat K_1$ form the
interval
\begin{equation}
\frac{a-\sigma\hat K_2}{\kappa}<\hat K_1<-\eta-(a+\kappa\eta)\hat K_2,
\label{eq:gain-interval}
\end{equation}
which is nonempty if and only if $D\hat K_2>a+\kappa\eta$.
\end{thm}

\begin{proof}
A real $2\times2$ matrix is Hurwitz if and only if its trace is negative and its
determinant positive, which are the two inequalities of \eqref{eq:gain-interval}. The
interval is nonempty if and only if
$a-\sigma\hat K_2<-\kappa\eta-\kappa(a+\kappa\eta)\hat K_2$, that is
$D\hat K_2>a+\kappa\eta$. If $D\neq0$ this holds for $\hat K_2$ of the sign of $D$ and
large magnitude. If $D=0$ then $a+\kappa\eta=\sigma/\kappa>0$ and the condition reads
$0>\sigma/\kappa$, which fails.
\end{proof}

\begin{rem}
Theorems~\ref{Possibility_theorem} and~\ref{Rule_theorem} together measure the cost of
deriving the rule from an objective. When $\rho\geq n$, an unrestricted rule can stabilize
whenever $D\neq0$, while a rule derived from a concave objective can do so only when
$D>0$. The difference is the sign restriction $BC>0$: by Lemma~\ref{Realization_lemma},
concavity bounds the third row of $J$ from below, and when $D<0$ every stabilizing third
row violates that bound. A static rule responding to consumption, $K_2\neq0$, plays the
role here that a nonzero $W_{mc}$ plays for the Central Bank: with $K_2=0$ the interval
\eqref{eq:gain-interval} is nonempty only if $a+\kappa\eta<0$, the static analogue of the
condition $d>\kappa\eta$ of Theorem~\ref{Possibility_kl}.
\end{rem}
\section{Conclusion}

We have added to the Ramsey--Cass--Koopmans system a Central Bank whose policy rule is
derived from its own welfare function, and studied the resulting three-dimensional system
under the assumption that both agents follow the differential form of their first-order
conditions from inherited initial values. The household system is controllable through
the Central Bank's instrument exactly when $D\neq0$
(Theorem~\ref{Controllability_theorem}). A Central Bank whose objective ignores
consumption cannot render the equilibrium locally asymptotically stable when $\rho\geq n$
(Theorem~\ref{Impossibility_theorem}); one whose objective responds to consumption can do
so if and only if $D>0$ (Theorem~\ref{Possibility_theorem}), which under constant relative
risk aversion and Cobb--Douglas technology is an upper bound on household risk aversion
(Proposition~\ref{Theta_bound}). When $\rho>n$, every trajectory converging to the stable
equilibrium satisfies both agents' transversality conditions, so that these conditions
no longer select the path (Theorem~\ref{TVC_theorem}). Oscillations arise only through
policy, and once $D>0$ cycles of every period and rate of decay are attainable
(Corollary~\ref{Cycles_iff}). 
Several questions remain open, with one of the perhaps more interesting ones being the feasibility of this framework in a multi-asset system proposed by Hahn \cite{Hahn1966}.


\begin{thebibliography}{99}

\bibitem{Barro1974}
R.~J. Barro,
\emph{Are Government Bonds Net Wealth?},
Journal of Political Economy \textbf{82} (1974), no.~6, 1095--1117.

\bibitem{Barro1990}
R.~J. Barro,
\emph{Government Spending in a Simple Model of Endogenous Growth},
Journal of Political Economy \textbf{98} (1990), no.~5, part 2, S103--S125.
\newblock \url{https://ideas.repec.org/a/ucp/jpolec/v98y1990i5ps103-26.html}

\bibitem{BenhabibNishimura1979}
J.~Benhabib and K.~Nishimura,
\emph{The Hopf Bifurcation and the Existence and Stability of Closed Orbits in Multisector Models of Optimal Economic Growth},
Journal of Economic Theory \textbf{21} (1979), no.~3, 421--444.

\bibitem{BrockScheinkman1976}
W.~A. Brock and J.~A. Scheinkman,
\emph{Global Asymptotic Stability of Optimal Control Systems with Applications to the Theory of Economic Growth},
Journal of Economic Theory \textbf{12} (1976), no.~1, 164--190.

\bibitem{Cass1965}
D.~Cass,
\emph{Optimum Growth in an Aggregative Model of Capital Accumulation},
Review of Economic Studies \textbf{32} (1965), no.~3, 233--240.

\bibitem{CoibionGorodnichenko2015}
O.~Coibion and Y.~Gorodnichenko,
\emph{Information Rigidity and the Expectations Formation Process: A Simple Framework and New Facts},
American Economic Review \textbf{105} (2015), no.~8, 2644--2678.

\bibitem{Domar1946}
E.~D. Domar,
\emph{Capital Expansion, Rate of Growth, and Employment},
Econometrica \textbf{14} (1946), no.~2, 137--147.

\bibitem{ElminejadHavranekIrsova2022}
A.~Elminejad, T.~Havranek and Z.~Irsova,
\emph{People Are Less Risk-Averse than Economists Think},
CEPR Discussion Paper No.~17411, Centre for Economic Policy Research, London, 2022.

\bibitem{EvansHonkapohja2001}
G.~W. Evans and S.~Honkapohja,
\emph{Learning and Expectations in Macroeconomics},
Princeton University Press, Princeton, NJ, 2001.
\newblock \url{https://www.jstor.org/stable/j.ctt7s6t9}

\bibitem{Hahn1960}
F.~H. Hahn,
\emph{The Stability of Growth Equilibrium},
Quarterly Journal of Economics \textbf{74} (1960), no.~2, 206--226.

\bibitem{Hahn1966}
F.~H. Hahn,
\emph{Equilibrium Dynamics with Heterogeneous Capital Goods},
Quarterly Journal of Economics \textbf{80} (1966), no.~4, 633--646.

\bibitem{HahnNegishi1962}
F.~H. Hahn and T.~Negishi,
\emph{A Theorem on Non-T\^atonnement Stability},
Econometrica \textbf{30} (1962), no.~3, 463--469.

\bibitem{Harrod1939}
R.~F. Harrod,
\emph{An Essay in Dynamic Theory},
Economic Journal \textbf{49} (1939), no.~193, 14--33.

\bibitem{Heller1975}
W.~P. Heller,
\emph{T\^atonnement Stability of Infinite Horizon Models with Saddle-Point Instability},
Econometrica \textbf{43} (1975), no.~1, 65--80.

\bibitem{Kamihigashi2001}
T.~Kamihigashi,
\emph{Necessity of Transversality Conditions for Infinite Horizon Problems},
Econometrica \textbf{69} (2001), no.~4, 995--1012.

\bibitem{Koopmans1965}
T.~C. Koopmans,
\emph{On the Concept of Optimal Economic Growth},
in \emph{The Econometric Approach to Development Planning},
Pontificiae Academiae Scientiarum Scripta Varia \textbf{28} (1965), 225--300.

\bibitem{Kurz1968}
M.~Kurz,
\emph{The General Instability of a Class of Competitive Growth Processes},
Review of Economic Studies \textbf{35} (1968), no.~2, 155--174.

\bibitem{LeiNoussair2002}
V.~Lei and C.~N. Noussair,
\emph{An Experimental Test of an Optimal Growth Model},
American Economic Review \textbf{92} (2002), no.~3, 549--570.
\newblock \url{https://www.aeaweb.org/articles?id=10.1257/00028280260136246}

\bibitem{Manski2004}
C.~F. Manski,
\emph{Measuring Expectations},
Econometrica \textbf{72} (2004), no.~5, 1329--1376.

\bibitem{Michel1982}
P.~Michel,
\emph{On the Transversality Condition in Infinite Horizon Optimal Problems},
Econometrica \textbf{50} (1982), no.~4, 975--986.

\bibitem{Ramsey1928}
F.~P. Ramsey,
\emph{A Mathematical Theory of Saving},
Economic Journal \textbf{38} (1928), no.~152, 543--559.

\bibitem{Solow1956}
R.~M. Solow,
\emph{A Contribution to the Theory of Economic Growth},
Quarterly Journal of Economics \textbf{70} (1956), no.~1, 65--94.

\bibitem{Tobin1969}
J.~Tobin,
\emph{A General Equilibrium Approach to Monetary Theory},
Journal of Money, Credit and Banking \textbf{1} (1969), no.~1, 15--29.

\bibitem{Turnovsky1997}
S.~J. Turnovsky,
\emph{Fiscal Policy in a Growing Economy with Public Capital},
Macroeconomic Dynamics \textbf{1} (1997), no.~3, 615--639.

\end{thebibliography}
\end{document}